\documentclass[12pt]{article}
\usepackage{amscd}
\usepackage{amsmath,amsfonts,amssymb,amscd}
\usepackage{indentfirst,graphicx,epsfig}
\usepackage{graphicx,psfrag}

\input{epsf}

\newtheorem{theo}{Theorem}
\newtheorem{coro}{Corollary}

\newtheorem{lemma}{Lemma}
\newtheorem{conj}{Conjecture}

\def\qed{\hfill \nopagebreak\rule{5pt}{8pt}}
\def\pf{\noindent {\it Proof.} }

\def\qed{\hfill \nopagebreak\rule{5pt}{8pt}}
\def\pf{\noindent {\it Proof.} }
\begin{document}

\title {\Large \bf Complete Solution to a Conjecture\\ on the Maximal Energy of Unicyclic
Graphs\footnote{Supported by NSFC and ``the Fundamental Research
Funds for the Central Universities". }}

\author{
  \small Bofeng Huo$^{1,2}$, Xueliang Li$^1$, Yongtang Shi$^1$ \\[3mm]
  \small $^1$Center for Combinatorics and LPMC-TJKLC\\
  \small Nankai University, Tianjin 300071, China\\
  \small E-mail: huobofeng@mail.nankai.edu.cn;  lxl@nankai.edu.cn; shi@nankai.edu.cn\\
  \small $^2$Department of Mathematics and Information Science\\
  \small Qinghai Normal University, Xining 810008, China
}

\date{}
\maketitle
\begin{abstract}

{\small For a given simple graph $G$, the energy of $G$, denoted by
$E(G)$, is defined as the sum of the absolute values of all
eigenvalues of its adjacency matrix. Let $P_n^{\ell}$ be the
unicyclic graph obtained by connecting a vertex of $C_\ell$ with a
leaf of $P_{n-\ell}$\,. In [G. Caporossi, D. Cvetkovi\'c, I. Gutman,
P. Hansen, Variable neighborhood search for extremal graphs. 2.
Finding graphs with extremal energy, {\it J. Chem. Inf. Comput.
Sci.\/} {\bf 39}(1999) 984--996], Caporossi et al. conjectured that
the unicyclic graph with maximal energy is $C_n$ if $n\leq 7$ and
$n=9,10,11,13,15$\,, and $P_n^6$ for all other values of $n$. In
this paper, by employing the Coulson integral formula and some
knowledge of real analysis, especially by using certain
combinatorial technique, we completely solve this conjecture.
However, it turns out that for $n=4$ the conjecture is not true, and
$P_4^3$ should be the unicyclic graph with maximal energy. \\
[2mm]
Keywords: graph energy; Coulson integral formula; unicyclic graph\\
AMS classification 2010: 05C50, 05C90, 15A18, 92E10.}

\end{abstract}

\section{Introduction}

For a given simple graph $G$ of order $n$, denote by $A(G)$ the
adjacency matrix of $G$. The characteristic polynomial of $A(G)$ is
usually called the characteristic polynomial of $G$, denoted by
\begin{equation*}\label{chapoly1}
\phi(G, x)=\det(xI-A(G))=x^n+a_1x^{n-1}+\cdots+a_n,
\end{equation*}
If $G$ is a bipartite graph, the characteristic polynomial of $G$
has the form
\begin{equation*}\label{cvetgraph}
\phi(G,
x)=\sum_{k=0}^{\lfloor{n/2}\rfloor}a_{2k}x^{n-2k}=\sum_{k=0}^{\lfloor{n/2}\rfloor}(-1)^{k}b_{2k}x^{n-2k},
\end{equation*}
where $b_{2k}=(-1)^ka_{2k}$ and $b_{2k}\geq 0$ for all
$k=1,\ldots,\lfloor{n/2}\rfloor$, especially $b_0=a_0=1$. In
particular, if $G$ is a tree, the characteristic polynomial of $G$
can be expressed as
\begin{equation*}\label{cvettree}
\phi(G,  x)=\sum_{k=0}^{\lfloor{n/ 2}\rfloor}(-1)^{k}m(G,k)x^{n-2k},
\end{equation*}
where $m(G,k)$ is the number of $k$-matchings of $G$.

For a graph $G$, let $\lambda_1, \lambda_2,\ldots, \lambda_n$ denote
the eigenvalues of $\phi(G, x)$. The {\it energy} of $G$ is defined
as
$$E(G)=\sum_{i=1}^n|\lambda_i|.$$
This definition was put forward by Gutman \cite{gutman1978} in 1978.
The following formula is also well-known
\begin{equation*}
E(G)={1\over\pi}\int^{+\infty}_{-\infty}{1\over x^2}\log |x^n
{\phi(G,i/x)}|\mathrm{d}x,
\end{equation*}
where $i^2=-1$. Furthermore, in the book of Gutman and Polansky
\cite{gutman&polansky1986}, the above equality was converted into an
explicit formula as follows:
\begin{equation*}\label{tujifen}
E(G)={1\over2\pi}\int^{+\infty}_{-\infty}{1\over
x^2}\log\left[\left(\sum_{k=0}^{\lfloor{n/
2}\rfloor}(-1)^ka_{2k}x^{2k}\right)^2+\left(\sum_{k=0}^{\lfloor{n/
2}\rfloor}(-1)^ka_{2k+1}x^{2k+1}\right)^2\right] \mathrm{d}x.
\end{equation*}
For more results about graph energy, we refer the readers to the
survey of Gutman, Li and Zhang \cite{gutman&lxl2009}.

For two given trees, or bipartite graphs $G_1$ and $G_2$, according
to the corresponding coefficients of the characteristic polynomials,
one can introduce a quasi order to compare the values of $E(G_1)$
and $E(G_2)$. Actually, the quasi order method is commonly used to
compare the energies of pairs of such graphs. However, for general
graphs, it is difficult to define such a quasi order. If, for two
trees, or bipartite graphs, the above quantities $m(T, k)$ or
$|a_k(G)|$ can not be compared uniformly, then the quasi order
method is invalid, and this happened very often. Recently, for these
quasi-order incomparable problems, we find an efficient approach to
determine which one attains the extremal value of the energy, such
as our earlier papers \cite{Huo&Ji&Li1}--\cite{Huo&Ji&Li3}.

Let $C_n$ be the cycle of order $n$, $P_n$ the path of order $n$,
and $P_n^{\ell}$ the unicyclic graph obtained by connecting a vertex
of $C_\ell$ with a leaf of $P_{n-\ell}$\,. In \cite{CC}, Caporossi
et al. proposed the following conjecture on the unicyclic graph with
maximal energy.
\begin{conj}\label{conjecture}
Among all unicyclic graphs on $n$ vertices, the cycle $C_n$ has
maximal energy if $n\leq 7$ and $n=9,10,11,13$ and $15$\,. For all
other values of $n$\,, the unicyclic graph with maximal energy is
$P_n^6$\,.
\end{conj}

In \cite{hou2002}, the authors proved the following Theorem
\ref{HLSthm0} that is weaker than the above conjecture, namely that
$E(P_n^6)$ is maximal within the class of the unicyclic bipartite
$n$-vertex graphs differing from $C_n$\,. And they also claimed that
the energy of $C_n$ and $P_n^6$ is quasi-order incomparable.
\begin{theo}\label{HLSthm0}
Let $G$ be any connected, unicyclic and bipartite graph on $n$
vertices and $G \ncong  C_n$\,. Then ${E}(G)<{E}(P_n^6)$\,.
\end{theo}

Very recently, our another paper \cite{HLS20101109} and Andriantiana
\cite{Andriantiana} independently proved that $E(C_n)<E(P_n^6)$, and
then completely determined that $P_n^6$ is the only graph which
attains the maximum value of the energy among all the unicyclic
bipartite graphs for $n=8,12,14$ and $n\geq 16$, which partially
solves the above conjecture.

\begin{theo}\label{HJLStheo}
For $n=8,12,14$ and $n\geq 16$, $E(P_n^6)>E(C_n)$.
\end{theo}

In this paper, by employing the Coulson integral formula (details on
the formula can be found in \cite{coulson} and
\cite{gutman&polansky1986} pp.139-147, as well as in the recent
works \cite{gut&mat, mat&gut}) and some knowledge of real analysis,
especially by using certain combinatorial technique, we completely
solve this conjecture by proving the following theorem and
corollary. However, we find that for $n=4$ the conjecture is not
true, and $P_4^3$ should be the unicyclic graph with maximal energy.

\begin{theo}\label{HLSthm20101109}
Among all unicyclic graphs of order $n\geq 16$\,, the unicyclic
graph with maximal energy is $P_n^6$\,.
\end{theo}

\begin{coro}\label{HLScoro20101109}
Among all unicyclic graphs on $n$ vertices, the cycle $C_n$ has
maximal energy if $n\leq 7$ but $n\neq 4$, and $n=9,10,11,13$ and
$15$\,;  $P_4^3$ has maximal energy if $n=4$\,. For all other values
of $n$\,, the unicyclic graph with maximal energy is $P_n^6$\,.
\end{coro}

\section{Preliminaries}

Let $G(n,\ell)$ be the set of all connected unicyclic graphs on $n$
vertices that contain the cycle $C_\ell$ as a subgraph. Denote by
$C(n,\ell)$ the set of all unicyclic graphs obtained from $C_\ell$
by adding to it $n-\ell$ pendent vertices. In the following, we list
some results given in \cite{hou2002} which will be used in the
sequel.

\begin{lemma}\label{thm20101109n}
Let $G\in G(n,\ell)$ and $n>\ell$\,. If $G$ has maximal energy in
$G(n,\ell)$\,, then $G$ is either $P_n^\ell$ or, when $\ell=4r$\,, a
graph from $C(n,\ell)$\,.
\end{lemma}

\begin{lemma}\label{LJlem3.7}
Let $G\in C(n,\ell)$ and $n>\ell$. If $\ell$ is even with $\ell\geq
8$ or $\ell=4$\,, then $E(G)< E(P_n^6)$\,.
\end{lemma}

\begin{lemma}\label{LJlem3.6}
Let $\ell$ be even and $\ell\geq 8$ or $\ell=4$\,. Then
$E(P_n^\ell)< E(P_n^6)$\,.
\end{lemma}

Form Lemmas \ref{thm20101109n}--\ref{LJlem3.6} and Theorem
\ref{HJLStheo}, we conclude that for any $n$-vertex unicyclic graph
$G$, if the length of the unique cycle of $G$ is even and
$n=8,12,14$ and $n\geq 16$, then $E(G)<E(P_n^6)$; if the length of
the unique cycle of $G$ is odd and $G\in G(n,\ell)$, then
$E(G)<E(P_n^\ell)$. For proving Theorem \ref{HLSthm20101109}, we
only need to show that $E(P_n^\ell)<E(P_n^6)$ for every odd $\ell$
and $n\geq 16$.

In the remainder of this section, we will introduce some lemmas and
notations. At first, we recall some knowledge on real analysis, for
which we refer the readers to \cite{zorich2002}.

\begin{lemma}\label{inequalitylemm}
For any real number $X>-1$, we have
\begin{equation*}
\frac{X}{1+X}\leq\log(1+ X)\leq X.
\end{equation*}
In particular, $\log(1+ X)<0$ if and only if $X<0$.
\end{lemma}
The following lemma on the difference of the energies of two graphs
is a well-known result due to Gutman \cite{gutman2001}, which will
be used in the sequel.
\begin{lemma}\label{lemma1.1}
If $G_1$ and $G_2$ are two graphs with the same number of vertices,
then
$$E(G_1)-E(G_2)={1\over\pi}\int^{+\infty}_{-\infty}
\log\left|\frac{ \phi (G_1, ix)}{\phi(G_2, ix)}\right|
\mathrm{d}x.$$
\end{lemma}

Now we present one basic formula of the characteristic polynomial
$\phi(G,x)$\,, which can be found in \cite{cvet1980}.

\begin{lemma} \label{zh8}
Let $uv$ be an edge of $G$\,. Then
$$
\phi (G,x)=\phi (G-uv,x) - \phi (G-u-v,x)-2\sum _{C\in {\mathcal
C}(uv)}\phi (G-C,x)
$$
where ${\mathcal C}(uv)$ is the set of cycles containing $uv$\,. In
particular, if $uv$ is a pendant edge with pendant vertex $v$\,,
then $ \phi (G,x)=x\,\phi (G-v,x)-\phi (G-u-v,x)$\,.
\end{lemma}
From Lemma \ref{zh8}, we can easily obtain the following lemma.
\begin{lemma}\label{HJLSlem1}
For any positive integer $t\leq n-2$,
$\phi(P_n^t,x)=x\phi(P_{n-1}^t,x)-\phi(P_{n-2}^t,x)$. In particular,
$\phi(P_n^6,x)=x\phi(P_{n-1}^6,x)-\phi(P_{n-2}^6,x)$.

\end{lemma}

Now for convenience, we introduce some notations as follows, which
will be well used in this sequel.
\begin{align*}
Y_1(x)=\frac{x+\sqrt{x^2-4}}{2},\qquad\qquad\qquad~
Y_2(x)=\frac{x-\sqrt{x^2-4}}{2}.
\end{align*}
It is easy to verify that $Y_1(x)+Y_2(x)=x$, $Y_1(x)Y_2(x)=1$,
$Y_1(ix)=\frac{x+\sqrt{x^2+4}}{2}i$ and
$Y_2(ix)=\frac{x-\sqrt{x^2+4}}{2}i$. We define
$$Z_1(x)=-iY_1(ix)=\frac{x+\sqrt{x^2+4}}{2},~
Z_2(x)=-iY_2(ix)=\frac{x-\sqrt{x^2+4}}{2}.$$ Observe that
$Z_1(x)+Z_2(x)=x$ and $Z_1(x)Z_2(x)=-1$. In addition, for $x>0$,
$Z_1(x)>1$ and $-1<Z_2(x)<0$; for $x<0$, $0<Z_1(x)<1$ and
$Z_2(x)<-1$. In the rest of this paper, we abbreviate $Z_j(x)$ to
$Z_j$ for $j=1,2$.

\section{Main results}

First, we introduce some more notations, which will be used
frequently later.
\begin{align*}
&A_1(x)=\frac {Y_1(x)\phi (P_8^6,x)-\phi (P_7^6,x)}
{(Y_1(x))^{9}-(Y_1(x))^{7}},\qquad A_2(x)=\frac {Y_2(x)\phi
(P_8^6,x)-\phi (P_7^6,x)} {(Y_2(x))^{9}-(Y_2(x))^{7}},\\
&B_1(x)=\frac {Y_1(x)\phi (P_{t+2}^t,x)-\phi (P_{t+1}^t,x)}
{(Y_1(x))^{t+3}-(Y_1(x))^{t+1}},~~ B_2(x)=\frac {Y_2(x)\phi
(P_{t+2}^t,x)-\phi (P_{t+1}^t,x)} {(Y_2(x))^{t+3}-(Y_2(x))^{t+1}},\\
&C_1(x)=\frac {Y_1(x)(x^2-1)-x}
{(Y_1(x))^{3}-Y_1(x)},\qquad\quad~~~~~~~ C_2(x)=\frac
{Y_2(x)(x^2-1)-x} {(Y_2(x))^{3}-Y_2(x)}.
\end{align*}
By some calculations, we can get that $\phi
(P_8^6,x)=x^8-8x^6+19x^4-16x^2+4$ and $\phi
(P_7^6,x)=x^7-7x^5+13x^3-7x$, and then
$$A_1(ix)=-\frac{Z_1f_8+f_7}{Z_1^2+1}Z_2^7, \quad
A_2(ix)=-\frac{Z_2f_8+f_7}{Z_2^2+1}Z_1^7,$$ where $f_8=\phi
(P_8^6,ix)=x^8+8x^6+19x^4+16x^2+4$ and $f_7=i\phi
(P_7^6,ix)=x^7+7x^5+13x^3+7x$.

\begin{lemma}\label{HJLSlem0}
For $n\geq 7$ and odd integer $3\leq t\leq n$, the characteristic
polynomials of $P_n^6$ and $P_n^t$ have the following forms:
$$\phi (P_n^6,x)=A_1(x)(Y_1(x))^n+A_2(x)(Y_2(x))^n$$ and
$$\phi (P_n^t,x)=B_1(x)(Y_1(x))^n+B_2(x)(Y_2(x))^n$$
where $x\neq \pm 2$.
\end{lemma}
\pf By Lemma \ref{HJLSlem1}, we notice that $\phi (P_n^6,x)$
satisfies the recursive formula $f(n,x)=xf(n-1,x)-f(n-2,x)$.
Therefore, the general solution of this linear homogeneous
recurrence relation is $f(n,x)=D_1(x)(Y_1(x))^n+D_2(x)(Y_2(x))^n$.
By some elementary calculations, we can easily obtain that
$D_i(x)=A_i(x)$ for $\phi (P_n^6,x)$, $i=1,2$, from the initial
values $\phi(P_8^6,x)$, $\phi(P_7^6,x)$. Similarly, the required
expression of $\phi (P_n^t,x)$ can be obtained by the analogous
method.\qed

Employing a method similar to the proof of Lemma \ref{HJLSlem0}, we
can obtain
\begin{lemma}\label{HLSlem1}
For positive integer $t\geq 3$, we have
\begin{align*}
\phi(P_{t+2}^t,x)=&\left(C_1(x)(Y_1(x))^{t-2}((Y_1(x))^4-x^2+1)\right)\\
&+\left(C_2(x)(Y_2(x))^{t-2}((Y_2(x))^4-x^2+1)\right)-2(x^2-1);\\[2mm]
\phi(P_{t+1}^t,x)=&\left(C_1(x)(Y_1(x))^{t-2}((Y_1(x))^3-x)\right)
+\left(C_2(x)(Y_2(x))^{t-2}((Y_2(x))^3-x)\right)-2x.
\end{align*}
\end{lemma}
\pf By Lemma \ref{zh8}, we notice that $\phi (P_n,x)$ satisfies the
recursive formula $f(n,x)=xf(n-1,x)-f(n-2,x)$. Therefore, the
general solution of this linear homogeneous recurrence relation is
$f(n,x)=D_1(x)(Y_1(x))^n+D_2(x)(Y_2(x))^n$. By some elementary
calculations, we can easily obtain that $D_i(x)=C_i(x)$ for $\phi
(P_n,x)$, $i=1,2$, from the initial values $\phi(P_2,x)$,
$\phi(P_1,x)$. According to Lemma \ref{zh8}, we have
\begin{align*}
\phi(P_{t+2}^t,x)=&\phi(P_{t+2},x)-\phi(P_{t-2},x)\phi(P_2,x)-2\phi(P_2,x);\\
\phi(P_{t+1}^t,x)=&\phi(P_{t+1},x)-\phi(P_{t-2},x)\phi(P_1,x)-2\phi(P_1,x).
\end{align*}
Therefore, we can obtain the required expression for
$\phi(P_{t+2}^t,x)$ and $\phi(P_{t+1}^t,x)$.\qed

Notice that $(x^2+1)Z_1+x=Z_1^3$ and $(x^2+1)Z_2+x=Z_2^3$. By some
simplifications, we can get the following corollary from Lemma
\ref{HLSlem1}.
\begin{coro}\label{HLScoro1}
$B_1(ix)=B_{11}(t,x)+B_{12}(t,x)\cdot i^t$ and
$B_2(ix)=B_{21}(t,x)+B_{22}(t,x)\cdot i^t$, where
\begin{align*}
&B_{11}(t,x)=\frac{Z_1^2(Z_1^2+2)}{(Z_1^2+1)^2}-\frac{Z_2^{2t-2}}{x^2+4},\qquad
B_{12}(t,x)=\frac{-2Z_2^{t-2}}{Z_1^2+1},\\[2mm]
&B_{21}(t,x)=\frac{Z_2^2(Z_2^2+2)}{(Z_2^2+1)^2}-\frac{Z_1^{2t-2}}{x^2+4},\qquad
B_{12}(t,x)=\frac{-2Z_1^{t-2}}{Z_2^2+1}.\\
\end{align*}
\end{coro}
For brevity of the exposition, we denote
$$g_1=\frac{Z_1^2(Z_1^2+2)}{(Z_1^2+1)^2},~~g_2=\frac{Z_2^2(Z_2^2+2)}{(Z_2^2+1)^2},~~
m_1=\frac{-2}{Z_1^2+1},~~m_2=\frac{-2}{Z_2^2+1},~~h=\frac{1}{x^2+4}.$$
Observe that each of $g_i$, $m_i$, $h$ is a real function only in
$x$, $i=1,2$.

From now on, we use $A_j$ and $B_{jk}$ instead of $A_j(ix)$ and
$B_{jk}(t,x)$ for $j,k=1, 2$, respectively. According to Lemma
\ref{HJLSlem0} and Corollary \ref{HLScoro1}, it is no hard to get
the following simplifications.
\begin{align}
&\left|\phi(P_n^6,ix)\right|^2=A_1^2Z_1^{2n}+A_{2}^2Z_2^{2n}+(-1)^n2A_1A_2\label{HLSeq3},\\[2mm]
&\left|\phi(P_n^t,ix)\right|^2=(B_{11}^2+B_{12}^2)Z_1^{2n}
+(B_{21}^2+B_{22}^2)Z_2^{2n}+(-1)^n2(B_{11}B_{21}+B_{12}B_{22}).\label{HLSeq2}
\end{align}

{\bf Proof of Theorem \ref{HLSthm20101109}.}

From the analysis in the above section, we only need to show that
$E(P_n^t)<E(P_n^6)$ for every odd $t\leq n$ and $n\geq 16$. By Lemma
\ref{lemma1.1},
$$E(P_n^t)-E(P_n^6)={1\over\pi}\int^{+\infty}_{-\infty}
\log\left|\frac{ \phi (P_n^t, ix)}{\phi(P_n^6, ix)}\right|
\mathrm{d}x.$$
We distinguish two cases in terms of the parity of
$n$.

{\bf Case 1.} $n$ is odd and $n\geq 17$.

Now we will prove that the integrand $\log\left|\frac{ \phi (P_n^t,
ix)}{\phi(P_n^6, ix)}\right|$ is monotonically decreasing in $n$.
\begin{align*}
&\log\left|\frac{ \phi (P_{n+2}^t, ix)}{\phi(P_{n+2}^6,
ix)}\right|-\log\left|\frac{ \phi (P_n^t, ix)}{\phi(P_n^6,
ix)}\right|\\[2mm]
=&\frac 1 2 \log\frac{\left| \phi (P_{n+2}^t, ix)\cdot\phi(P_n^6,
ix)\right|^2}{\left|\phi(P_{n+2}^6, ix)\cdot\phi (P_n^t,
ix)\right|^2} =\frac 1 2\log\left(1+\frac {K(n,t,x)}
{H(n,t,x)}\right),
\end{align*}
where $H(n,t,x)=\left|\phi(P_{n+2}^6, ix)\cdot\phi (P_n^t,
ix)\right|^2>0$ and $$K(n,t,x)=\left| \phi (P_{n+2}^t,
ix)\cdot\phi(P_n^6, ix)\right|^2-\left|\phi(P_n^6, ix)\cdot\phi
(P_n^t, ix)\right|^2.$$ From Lemma \ref{inequalitylemm}, we only
need to prove $K(n,t,x)<0$. By some elementary calculations and
simplifications, we can obtain
\begin{align*}
K(n,t,x)=\alpha(t,x)(Z_1^4-Z_2^4)+\beta(t,x)Z_1^{2n}(Z_1^{4}-1)+\gamma(t,x)Z_2^{2n}(1-Z_2^{4}),
\end{align*}
where
$\alpha(t,x)=A_2^2(B_{11}^2+B_{12}^2)-A_1^2(B_{21}^2+B_{22}^2)$,
$\beta(t,x)=2A_1^2(B_{11}B_{21}+B_{12}B_{22})-2A_1A_2(B_{11}^2+B_{12}^2)$,
$\gamma(t,x)=2A_1A_2(B_{21}^2+B_{22}^2)-2A_2^2(B_{11}B_{21}+B_{12}B_{22})$.
In the following, we will discuss the signs of
$\alpha(t,x),\,\beta(t,x),\,\gamma(t,x)$.
\begin{align*}
\alpha(t,x)=&~\alpha_0+\alpha_1Z_1^{2t-4}+\alpha_2Z_2^{2t-4}+\alpha_3Z_1^{4t-4}+\alpha_4Z_2^{4t-4},\\[2mm]
\beta(t,x)=&~\beta_0+\beta_1Z_1^{2t-2}+\beta_2Z_2^{2t-2}\qquad\quad\quad~~+\beta_4Z_2^{4t-4},\\[2mm]
\gamma(t,x)=&~\gamma_0+\gamma_1Z_1^{2t-2}+\gamma_2Z_2^{2t-2}+\gamma_3Z_1^{4t-4},
\end{align*}
where
\begin{align*}
&\alpha_0=A_2^2g_1^2-A_1^2g_2^2,\qquad\qquad\qquad\quad\quad~~~\alpha_1=2A_1^2g_2hZ_1^2-A_1^2m_2^2,\\[2mm]
&\alpha_2=A_2^2m_1^2-2A_2^2g_1hZ_2^2,~\quad\qquad\quad\quad~~~\alpha_3=-A_1^2h^2,~\qquad\alpha_4=A_2^2h^2,\\[2mm]
&\beta_0=-2A_1\left(\frac{2(x^2+3)}{(x^2+4)^2}A_1+A_2g_1^2\right),\quad~\beta_1=-2A_1^2g_1h,\\[2mm]
&\beta_2=2A_1(2A_2g_1h-A_1g_2h-A_2m_1^2Z_1^2),~\beta_4=-2A_1A_2h^2,\\[2mm]
&\gamma_0=2A_2\left(A_1g_2^2+\frac{2(x^2+3)}{(x^2+4)^2}A_2\right),~~~~\quad\gamma_1=2A_2(A_1m_2^2Z_2^2+A_2g_1h-2A_1g_2h),\\[2mm]
&\gamma_2=2A_2^2g_2h\,,\qquad\quad\qquad\quad\qquad\quad\quad\quad~\gamma_3=2A_1A_2h^2.\\[2mm]
\end{align*}
{\bf Claim 1.} For any real number $x$ and positive integer $t$,
$\beta(t,x)<0$.

Notice that $Z_1f_8+f_7=(\frac x 2 f_8+f_7)+\frac
{\sqrt{x^2+4}}{2}f_8$, $Z_2f_8+f_7=(\frac x 2 f_8+f_7)-\frac
{\sqrt{x^2+4}}{2}f_8$ and
$$\left(\frac x 2 f_8+f_7\right)^2-\left(\frac
{\sqrt{x^2+4}}{2}f_8\right)^2=-(x^{10}+10x^8+36x^6+62x^4+51x^2+16)<0.$$
Then $A_1=-\frac{Z_1f_8+f_7}{Z_1^2+1}Z_2^7>0$,
$A_2=-\frac{Z_2f_8+f_7}{Z_2^2+1}Z_1^7>0$ since $Z_1>0$ and $Z_2<0$.
Therefore, $\beta_0<0$.
\begin{align*}
\beta_2=-\frac{A_1(x^2+1)}{(x^2+4)^{\frac 5
2}}&(x^9+11x^7+47x^5+93x^3+74x\\
&+\sqrt{x^2+4}(3x^8+27x^6+85x^4+111x^2+52))<0,
\end{align*}
since
\begin{align}
(x^9+11x^7+47x^5+93x^3+74x)^2-
(x^2+4)(3x^8+27x^6+85x^4+111x^2+52)^2<0.\label{HLSeq1}
\end{align}
It is easy to check that $\beta_1<0$ and $\beta_4<0$. Hence, the
claim holds.

{\bf Claim 2.} For any real number $x$ and positive integer $t$,
$\gamma(t,x)>0$.

Analogously, we can get $\gamma_0>0$, $\gamma_2>0$ and $\gamma_3>0$.
From Eq. \eqref{HLSeq1}, we have
\begin{align*}
\gamma_1=\frac{A_2(x^2+1)}{(x^2+4)^{\frac 5
2}}&(-(x^9+11x^7+47x^5+93x^3+74x)\\
&+\sqrt{x^2+4}(3x^8+27x^6+85x^4+111x^2+52))>0.
\end{align*}
Therefore, $\gamma(t,x)>0$.

{\bf Claim 3.} For any real number $x$ and odd $n\geq t$,
$K(n,t,x)\leq
\alpha(t,x)(Z_1^4-Z_2^4)+\beta(t,x)Z_1^{2t}(Z_1^{4}-1)+\gamma(t,x)Z_2^{2t}(1-Z_2^{4})$.

Since $Z_1(x)>1$ and $-1<Z_2(x)<0$ for $x>0$, we have $Z_1^{2n}\geq
Z_1^{2t}$ and $Z_2^{2n}\leq Z_2^{2t}$ when $n\geq t$. Since
$0<Z_1(x)<1$ and $Z_2(x)<-1$ for $x<0$, we have $Z_1^{2n}\leq
Z_1^{2t}$ and $Z_2^{2n}\geq Z_2^{2t}$ when $n\geq t$. From Claims 1
and 2, we have $\beta(t,x)<0$ and $\gamma(t,x)>0$ for any real
number $x$. Thus, Claim 3 holds.

{\bf Claim 4.}
$f(t,x)=\alpha(t,x)(Z_1^4-Z_2^4)+\beta(t,x)Z_1^{2t}(Z_1^{4}-1)+\gamma(t,x)Z_2^{2t}(1-Z_2^{4})$
is monotonically decreasing in $t$.

It is no difficult to get that
$f(t,x)=d_0+d_1Z_1^{2t}+d_2Z_2^{2t}+d_3Z_1^{4t}+d_4Z_2^{4t}=d_0+d_1(Z_1^2)^{t}+d_2(Z_1^2)^{-t}+d_3(Z_1^2)^{2t}
+d_4(Z_1^2)^{-2t}$, where
\begin{align*}
&d_0=\alpha_0(Z_1^4-Z_2^4)+\beta_2(Z_1^4-1)Z_1^2+\gamma_1(1-Z_2^4)Z_2^2,\\[2mm]
&d_1=\alpha_1(1-Z_2^8)+\beta_0(Z_1^4-1)+\gamma_3(Z_2^4-Z_2^8),\\[2mm]
&d_2=\alpha_2(Z_1^8-1)+\gamma_0(1-Z_2^4)+\beta_4(Z_1^8-Z_1^4),\\[2mm]
&d_3=\alpha_3(1-Z_2^8)+\beta_1(Z_1^2-Z_2^2),\\[2mm]
&d_4=\alpha_4(Z_1^8-1)+\gamma_2(Z_1^2-Z_2^2).
\end{align*}
We define $p_1(x)=x^3+6x$, $q_1(x)=(3x^2+4)\sqrt{x^2+4}$,
$p_2(x)=x^7+9x^5+24x^3+18x$,
$q_2(x)=(x^6+7x^4+12x^2+4)\sqrt{x^2+4}$,
$p_3(x)=x^{13}+15x^{11}+89x^9+264x^7+405x^5+288x^3+56x$,
$q_3(x)=(x^{12}+15x^{10}+85x^8+234x^6+331x^4+220x^2+48)\sqrt{x^2+4}$.
By some calculations, we have
\begin{align*}
&d_1=\frac{x(x^2+4)(x^2+1)^2(x-\sqrt{x^2+4})^7(p_2(x)+q_2(x))(p_3(x)+q_3(x))}{4(x^2+4-x\sqrt{x^2+4})^2(x^2+4+x\sqrt{x^2+4})^4},\\[2mm]
&d_2=\frac{x(x^2+4)(x^2+1)^2(x+\sqrt{x^2+4})^7(p_2(x)-q_2(x))(p_3(x)-q_3(x))}{4(x^2+4+x\sqrt{x^2+4})^2(x^2+4-x\sqrt{x^2+4})^4},\\[2mm]
&d_3=-\frac{x(x^2+1)^2(x-\sqrt{x^2+4})^{14}(p_1(x)+q_1(x))(p_2(x)+q_2(x))^2}{8192(x^2+4+x\sqrt{x^2+4})^4},\\[2mm]
&d_4=-\frac{x(x^2+1)^2(x+\sqrt{x^2+4})^{14}(p_1(x)-q_1(x))(p_2(x)-q_2(x))^2}{8192(x^2+4-x\sqrt{x^2+4})^4}.
\end{align*}
Since $(p_1(x))^2-(q_1(x))^2<0$, $(p_2(x))^2-(q_2(x))^2<0$ and
$(p_3(x))^2-(q_3(x))^2<0$, we deduce that, $d_1,\,d_3<0$ and
$d_2,\,d_4>0$ for $x>0$; $d_1,\,d_3>0$ and $d_2,\,d_4<0$ for $x<0$.
Therefore, no matter what of $x>0$ or $x<0$ happens, we always have
$$\frac {\partial f(t,x)}{\partial t}=\left(d_1(Z_1^2)^{t}-d_2(Z_1^2)^{-t}+2d_3(Z_1^2)^{2t}
-2d_4(Z_1^2)^{-2t}\right)\log Z_1^2<0.$$ The proof of Claim 4 is
complete.

From Claim 4, it follows that for $t\geq 5$, we have
\begin{align*}K(n,t,x)\leq f(5,x)
=&~-{x^2(x^2+1)^2(x^4+3x^2+1)}\\[2mm]
&\cdot (2x^{12}+31x^{10}+189x^8+574x^6+899x^4+661x^2+160) <0.
\end{align*}
For $t=3$, one must have $n>t+2$. So
\begin{align*}
&~K(n,3,x)<\alpha(3,x)(Z_1^4-Z_2^4)+\beta(3,x)Z_1^{2\times
3+4}(Z_1^{4}-1)+\gamma(3,x)Z_2^{2\times
3+4}(1-Z_2^{4})\\[2mm]
=&~-{x^2(x^2+1)^3(x^2+5)}
(2x^{12}+23x^{10}+104x^8+238x^6+290x^4+171x^2+32) <0.
\end{align*}
\begin{table}[ht]
\begin{center}
\arrayrulewidth 0.5 pt
\begin{tabular}{|p{30pt}|p{100pt}||p{30pt}|p{100pt}|}
\hline
 \quad $t$ & \quad $E(P_{17}^t)-E(P_{17}^6)$ &~~~ $t$ & \quad $E(P_{17}^t)-E(P_{17}^6)$  \\
\hline   \quad $3$  &\qquad $-0.05339$  & \quad $11$&\qquad $-0.12030$  \\
\hline   \quad $5$  &\qquad $-0.09835$ &\quad $13$ &\qquad  $-0.11425$ \\
\hline   \quad $7$  &\qquad $-0.11405$  &\quad $15$&\qquad  $-0.09493$ \\
\hline   \quad $9$  &\qquad $-0.12006$  & &   \\
\hline
\end{tabular}
\end{center}
\caption{The values of $E(P_{17}^t)-E(P_{17}^6)$ for $t\leq 15$.}
\label{HLS20101114table1}
\end{table}
We conclude that the integrand $\log\left|\frac{ \phi (P_n^t,
ix)}{\phi(P_n^6, ix)}\right|$ is monotonically decreasing in $n$.
Therefore, by Theorem \ref{HJLStheo}, for $n\geq 17$ and $t\geq 17$,
$E(P_n^t)-E(P_n^6)<E(P_{t}^t)-E(P_{t}^6)<0$. For $n\geq 17$ and
$t\leq 15$, $E(P_n^t)-E(P_n^6)<E(P_{17}^t)-E(P_{17}^6)<0$ from Table
\ref{HLS20101114table1}.

{\bf Case 2.} $n$ is even and $n\geq 8$.

From Eqs. \eqref{HLSeq2} and \eqref{HLSeq3}, we have
$$
\log\left|\frac{ \phi (P_n^t, ix)}{\phi(P_n^6, ix)}\right|^2=  \log
\frac{(B_{11}^2+B_{12}^2)Z_1^{2n}
+(B_{21}^2+B_{22}^2)Z_2^{2n}+2(B_{11}B_{21}+B_{12}B_{22})}{A_1^2Z_1^{2n}+A_{2}^2Z_2^{2n}+2A_1A_2}.$$
Therefore, when $n\rightarrow \infty$, we have
$$\left|\frac{ \phi (P_n^t, ix)}{\phi(P_n^6, ix)}\right|^2\rightarrow\left\{
 \begin{array}{ll}
 \frac{B_{11}^2+B_{12}^2}{A_1^2} &\mbox{if $x>0$}\\[3mm]
\frac{B_{21}^2+B_{22}^2}{A_2^2} &\mbox{if $x<0$}.
 \end{array}
 \right.
 $$
In this case, we will show $$\log \left|\frac{ \phi (P_n^t,
ix)}{\phi(P_n^6, ix)}\right|^2<\log
\frac{B_{11}^2+B_{12}^2}{A_1^2}$$ for $x>0$, and  $$\log
\left|\frac{ \phi (P_n^t, ix)}{\phi(P_n^6, ix)}\right|^2<\log
\frac{B_{21}^2+B_{22}^2}{A_2^2} $$ for $x<0$. Now we can simplify
the expressions of $\alpha_i$ for $i=0,1,2$ as follows:
\begin{align*}
&\alpha_0=\frac{x(x^2+1)^2(x^8+11x^6+43x^4+73x^2+50)(x^8+9x^6+27x^4+33x^2+12)}{(x^2+4)^{5/2}},\\[2mm]
&\alpha_1=-\frac{(p_2(x)+q_2(x))^2(3x^2+10+x\sqrt{x^2+4})(x-\sqrt{x^2+4})^{14}(x^2+1)^2}
{4096(x^2-x\sqrt{x^2+4}+4)^2(x^2+x\sqrt{x^2+4}+4)^2(x^2+4)},\\[2mm]
&\alpha_2=\frac{(p_2(x)-q_2(x))^2(3x^2+10-x\sqrt{x^2+4})(x+\sqrt{x^2+4})^{14}(x^2+1)^2}
{4096(x^2-x\sqrt{x^2+4}+4)^2(x^2+x\sqrt{x^2+4}+4)^2(x^2+4)}.
\end{align*}

{\bf Subcase 2.1.} $x>0$.

By some calculations, we have
$$\log \left|\frac{ \phi (P_n^t,
ix)}{\phi(P_n^6, ix)}\right|^2-\log
\frac{B_{11}^2+B_{12}^2}{A_1^2}=\log \left(1+\frac {K_1(n,t,x)}
{H_1(n,t,x)}\right),$$ where $H_1(n,t,x)=\left|\phi(P_n^6,
ix)\right|^2(B_{11}^2+B_{12}^2)>0$ and
$K_1(n,t,x)=-\alpha(t,x)Z_2^{2n}+\beta(t,x)$. Now we suppose
$\alpha(t,x)<0$. Otherwise, $K_1(n,t,x)<0$ since $\beta(t,x)<0$ by
Claim 1, and then we are done. Since $-1<Z_2<0$,
\begin{align*}
K_1(n,t,x)\leq
-\alpha(t,x)Z_2^{2t}+\beta(t,x)=\overline{d}_0+\overline{d}_1Z_1^{2t-2}+\overline{d}_2Z_2^{2t-2}
+\overline{d}_3Z_2^{4t-4}+\overline{d}_4Z_2^{6t-4},
\end{align*}
where $\overline{d}_0=\beta_0-\alpha_1Z_2^4$,
$\overline{d}_1=\beta_1-\alpha_3Z_2^2$,
$\overline{d}_2=\beta_2-\alpha_0Z_2^2$,
$\overline{d}_3=\beta_4-\alpha_2$, $\overline{d}_4=-\alpha_4$. Since
$\beta_i<0$ for $i=0,1,2,4$, $\alpha_0,\,\alpha_2,\,\alpha_4>0$ and
$\alpha_1,\,\alpha_3<0$, we have $\overline{d}_i<0$ for $i=2,3,4$
and
$$\overline{d}_1=-2A_1^2g_1h+A_1^2h^2Z_2^2=A_1^2h(hZ_2^2-2g_1)=-\frac{A_1^2h(2Z_1^2-Z_2^2+4)}{x^2+4}<0.$$
Denote by
$p_0(x)=x^{14}+19x^{12}+146x^{10}+584x^8+1300x^6+1582x^4+928x^2+160$
and
$q_0(x)=(x^{13}+17x^{11}+116x^9+404x^7+756x^5+722x^3+272x)\sqrt{x^2+4}$.
Then,
$$\overline{d}_0=-\frac{A_1(x^2+1)}{(Z_1^2+1)^4(Z_2^2+1)^2}\left(p_0(x)+q_0(x)\right)<0.$$
Thus, for $x>0$, $K_1(n,t,x)<0$, and
then $$\log \left|\frac{ \phi (P_n^t, ix)}{\phi(P_n^6,
ix)}\right|^2<\log \frac{B_{11}^2+B_{12}^2}{A_1^2}.$$

{\bf Subcase 2.2.} $x<0$.

Similarly, we can obtain
$$\log \left|\frac{ \phi (P_n^t,
ix)}{\phi(P_n^6, ix)}\right|^2-\log
\frac{B_{21}^2+B_{22}^2}{A_2^2}=\log \left(1+\frac {K_2(n,t,x)}
{H_2(n,t,x)}\right),$$ where $H_2(n,t,x)=\left|\phi(P_n^6,
ix)\right|^2(B_{21}^2+B_{22}^2)>0$ and
$K_2(n,t,x)=\alpha(t,x)Z_1^{2n}-\gamma(t,x)$. Now we suppose
$\alpha(t,x)>0$. Otherwise, $K_2(n,t,x)<0$ since $\gamma(t,x)>0$ by
Claim 2, and then we are done. Since $0<Z_1<1$,
\begin{align*}
K_2(n,t,x)\leq
\alpha(t,x)Z_1^{2t}-\gamma(t,x)=\widetilde{d}_0+\widetilde{d}_1Z_1^{2t-2}+\widetilde{d}_2Z_2^{2t-2}
+\widetilde{d}_3Z_1^{4t-4}+\widetilde{d}_4Z_1^{6t-4},
\end{align*}
where $\widetilde{d}_0=\alpha_2Z_1^4-\gamma_0$,
$\widetilde{d}_1=\alpha_0Z_1^2-\gamma_1$,
$\widetilde{d}_2=\alpha_4Z_1^2-\gamma_2$,
$\widetilde{d}_3=\alpha_1-\gamma_3$, $\widetilde{d}_4=\alpha_3$.
Since $\gamma_i>0$ for $i=0,1,2,3$,
$\alpha_0,\,\alpha_1,\,\alpha_3<0$ and $\alpha_2,\,\alpha_4>0$, we
have  $\widetilde{d}_i<0$ for $i=1,3,4$ and
$$\widetilde{d}_0=-\frac{A_2(x^2+1)}{(Z_2^2+1)^4(Z_1^2+1)^2}\left(p_0(x)-q_0(x)\right)<0,$$
$$\widetilde{d}_2=A_2^2h^2Z_1^2-2A_2^2g_2h=A_2^2h(hZ_1^2-2g_2)=-\frac{A_2^2h(2Z_2^2-Z_1^2+4)}{x^2+4}<0.$$
Thus, for $x<0$, $K_2(n,t,x)<0$, and then $$\log \left|\frac{ \phi
(P_n^t, ix)}{\phi(P_n^6, ix)}\right|^2<\log
\frac{B_{21}^2+B_{22}^2}{A_2^2}.
$$

From the two subcases, we conclude that
\begin{align*}
E(P_n^t)-E(P_n^6)&~={1\over\pi}\int^{+\infty}_{-\infty}
\log\left|\frac{ \phi (P_n^t, ix)}{\phi(P_n^6, ix)}\right|
\mathrm{d}x\\[2mm]
&~={1\over 2\pi}\int^{+\infty}_{-\infty}\log \left|\frac{
\phi (P_n^t, ix)}{\phi(P_n^6, ix)}\right|^2\mathrm{d}x\\[2mm]
&~<{1\over 2\pi}\int^{+\infty}_{0}\log
\frac{B_{11}^2+B_{12}^2}{A_1^2}\mathrm{d}x+{1\over
2\pi}\int^{0}_{-\infty}\log
\frac{B_{21}^2+B_{22}^2}{A_2^2}\mathrm{d}x.
\end{align*}
Denote
$p_4(x)=x^{16}+14x^{14}+83x^{12}+274x^{10}+551x^8+686x^6+507x^4+190x^2+22$,
$q_4(x)=(x^{15}+12x^{13}+61x^{11}+172x^9+291x^7+296x^5+167x^3+40x)\sqrt{x^2+4}$.
Notice that $\frac {Z_1^2}{(Z_1^2+1)^2}=\frac
{Z_2^2}{(Z_2^2+1)^2}=\frac 1 {x^2+4}$ and
$(p_4(x))^2-(q_4(x))^2=4(x^2+1)^2(2x^{10}+24x^8+104x^6+225x^4+248x^2+121)>0$
whenever $x>0$ or $x<0$. When $x>0$, $Z_2^2<1$, we have
\begin{align*}
B_{11}^2+B_{12}^2-A_1^2&=\left(\frac{Z_1^2+2}{x^2+4}-\frac{Z_2^{2t-2}}{x^2+4}\right)^2
+\left(-\frac{2(Z_1^2+1)Z_2^t}{x^2+4}\right)^2-A_1^2\\[2mm]
&=\frac{1}{(x^2+4)^2}\left((Z_1^2+2)^2+(2Z_1^2+4Z_2^2+4)Z_2^{2t-2}+Z_2^{4t-4}\right)-A_1^2\\
&<\frac{1}{(x^2+4)^2}\left((Z_1^2+2)^2+(2Z_1^2+4Z_2^2+4)Z_2^4+Z_2^8\right)-A_1^2\\
&=-\frac{p_4(x)-q_4(x)}{(x^2+4)(x^2+2+x\sqrt{x^2+4})}<0.
\end{align*}
When $x<0$, $Z_1^2<1$, we have
\begin{align*}
B_{21}^2+B_{22}^2-A_2^2&=\left(\frac{Z_2^2+2}{x^2+4}-\frac{Z_1^{2t-2}}{x^2+4}\right)^2
+\left(-\frac{2(Z_2^2+1)Z_1^t}{x^2+4}\right)^2-A_2^2\\[2mm]
&=\frac{1}{(x^2+4)^2}\left((Z_2^2+2)^2+(2Z_2^2+4Z_1^2+4)Z_1^{2t-2}+Z_1^{4t-4}\right)-A_2^2\\
&<\frac{1}{(x^2+4)^2}\left((Z_2^2+2)^2+(2Z_2^2+4Z_1^2+4)Z_1^4+Z_1^8\right)-A_2^2\\
&=-\frac{p_4(x)+q_4(x)}{(x^2+4)(x^2+2-x\sqrt{x^2+4})}<0.
\end{align*}
So
$$\int^{+\infty}_{0}\log
\frac{B_{11}^2+B_{12}^2}{A_1^2}\mathrm{d}x<0~~\mathrm{and}~~\int^{0}_{-\infty}\log
\frac{B_{21}^2+B_{22}^2}{A_2^2}\mathrm{d}x<0.$$
\begin{figure}[ht]
\begin{center}
\includegraphics[width=12cm]{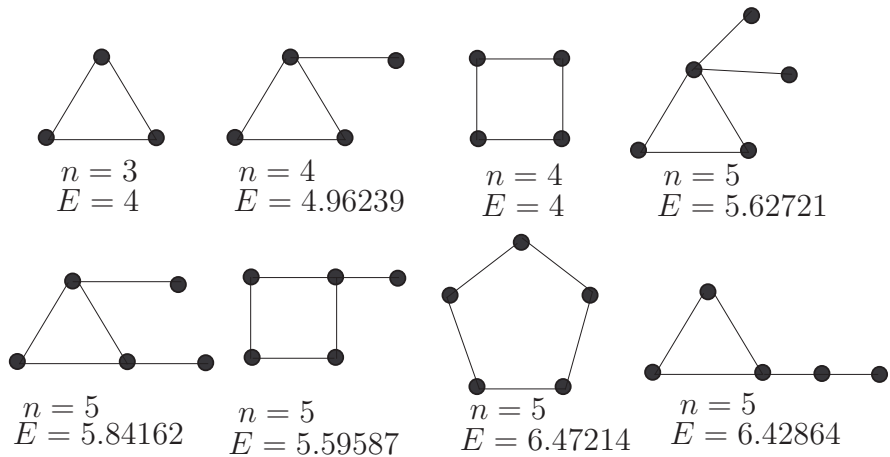}
\end{center}
\caption{All unicyclic graphs and its energies for $n\leq 5$.}
\label{fig}
\end{figure}

Therefore, $E(P_n^t)-E(P_n^6)<0$ when $n$ is even. \qed

\begin{table}[ht]
\begin{center}
\arrayrulewidth 0.5 pt
\begin{tabular}{|p{30pt}|p{30pt}|p{100pt}||p{30pt}|p{30pt}|p{100pt}|}
\hline
 \quad $n$ & \quad $t$ & \quad $E(P_{n}^t)-E(P_{n}^6)$ &~~  $n$& \quad $t$ & \quad $E(P_{n}^t)-E(P_{n}^6)$  \\
\hline \quad $6$  &\quad $3$  &\quad\quad $-0.45075 $  & \quad$6 $ &\quad $5$&\quad\quad $-0.53412 $  \\
\hline \quad $7 $  &\quad $3$  &\qquad\quad $0.22026 $ & \quad$7 $&\quad $5$ &\qquad \quad $0.19680$ \\
\hline \quad  $8 $ &\quad $3$  &\quad\quad $ -0.31283$  & \quad$8 $&\quad $5$&\quad\quad   $-0.37252$ \\
\hline \quad $ 8$  &\quad $7$  &\quad\quad $-0.42994 $  & \quad$9 $&\quad $3$&\qquad \quad $0.08604$ \\
\hline \quad  $ 9$ &\quad $5$  &\qquad\quad $0.04987 $  & \quad$9 $&\quad $7$&\qquad \quad $0.05443$ \\
\hline ~ $10 $ &\quad $3$  &\quad\quad $-0.26573 $  & ~ $10 $&\quad $5$&\quad\quad $-0.31918$ \\
\hline ~  $ 10$ &\quad $7$  &\quad\quad $-0.35115 $ & ~ $10 $&\quad $9$&\quad \quad $-0.40167$ \\
\hline ~  $11 $ &\quad $3$  &\qquad\quad $0.02396 $  & ~ $11 $&\quad $5$&\quad\quad $-0.01682$ \\
\hline ~  $11 $ &\quad $7$  &\quad\quad $-0.02469$  & ~ $11 $&\quad $9$&\quad \quad $-0.01186$ \\
\hline ~ $12$ &\quad $3$   &\quad\quad $-0.24081$  & ~ $12 $&\quad $5$&\quad\quad  $-0.29174$ \\
\hline ~  $12$ &\quad $7$   &\quad\quad $-0.31698$  & ~ $ 12$&\quad $9$&\quad \quad $-0.34102$ \\
\hline~  $12$ &~\,  $11$  &\quad\quad $-0.38894$  & ~ $13 $&\quad $3$&\quad\quad  $-0.01237$ \\
\hline ~  $13$ &\quad $5$   &\quad\quad $-0.05536$  & ~ $13 $&\quad $7$&\quad \quad $-0.06773 $ \\
\hline ~  $13$ &\quad $9$  &\quad\quad $-0.06719$  & ~ $13 $&~\,  $11$&\quad\quad  $-0.05081$ \\
\hline ~  $14$ &\quad $3$  &\quad \quad$-0.22520 $  & ~ $14 $&\quad $5$&\quad\quad  $ -0.27486$ \\
\hline ~  $14$ &\quad $7$  &\quad\quad $-0.29740 $  & ~ $ 14$&\quad $9$&\quad \quad $-0.31438$ \\
\hline ~  $14$ &~\, $11$  &\quad\quad $-0.33517$  & ~ $ 14$&~\,  $13$&\quad \quad $-0.38193$ \\
\hline ~  $15$ &\quad $3$  &\quad\quad $-0.03635$  & ~ $ 15$&\quad $5$&\quad\quad  $ -0.08055$ \\
\hline ~  $15$ &\quad $7$  &\quad\quad $-0.09506$  & ~ $ 15$&\quad $9$&\quad \quad $-0.09897$ \\
\hline ~  $15$ &~\, $11$  &\quad\quad $-0.09481 $  & ~ $15 $&~\,  $13$&\quad \quad $-0.07658$ \\
\hline ~  $16$ &\quad $3$  &\quad\quad $-0.21447$  & ~ $ 16$&\quad $5$&\quad \quad $-0.26340 $ \\
\hline ~  $16$ &\quad $7$  &\quad \quad$-0.28459 $  & ~ $16 $&\quad $9$&\quad \quad $-0.29873$ \\
\hline ~  $16$ &~\, $11$  &\quad\quad $-0.31223$  & ~ $16 $&~\,  $13$&\quad \quad $-0.33141$ \\
\hline ~  $16$ &~\, $15$  &\quad\quad $-0.37761$  & ~ $  $&~\,  $ $& \\
\hline
\end{tabular}
\end{center}
\caption{Values of $E(P_n^t)-E(P_n^6)$ for $n\leq 16$ and odd $t$.}
\label{HLS20101114table2}
\end{table}

{\bf Proof of Corollary \ref{HLScoro20101109}.}

There are only two unicyclic graphs of order $4$, which are shown in
Figure \ref{fig}. Observe that $P_4^3$ has maximal energy for $n=4$.
From Lemmas \ref{thm20101109n}--\ref{LJlem3.6}, and Theorems
\ref{HJLStheo} and \ref{HLSthm20101109}, we only need to show that
for $n\leq 16$ $(n\neq 4)$ and any odd $t$ with $3\leq t\leq n$\,,
$E(P_n^t)<E(P_n^6)$ or $E(P_n^t)<E(C_n)$\,. From Table
\ref{HLS20101114table2}, we can see that $E(P_n^t)<E(P_n^6)$ for
$6\leq n\leq 16$ except for $n=7,9,11$ and some $t$. In such cases,
we can check that $E(P_n^t)<E(C_n)$ from Table
\ref{HLS20101114table3}. For $n=3, 5$, we consider all the unicyclic
graphs. All such graphs and their energies are shown in Figure
\ref{fig}, in which our results are verified. Finally, we calculate
the energies of $C_n$ and $P_n^6$ for $n=7,9,10,11,13,15$, and
verify that $E(C_n)>E(P_n^6)$ in these cases.\qed

\begin{table}[ht]
\begin{center}
\arrayrulewidth 0.5 pt
\begin{tabular}{|p{30pt}|p{30pt}|p{40pt}|p{40pt}||p{30pt}|p{30pt}|p{40pt}|p{40pt}|}
\hline
\quad $n$ & \quad $t$ &  ~$E(P_{n}^t)$ &  ~$E(C_{n})$ &~~  $n$& \quad $t$ & ~$E(P_{n}^t)$ & ~$E(C_{n})$ \\
\hline \quad $7 $  &\quad $3$  & $8.94083$ &  $8.98792$ & \quad$7$ &\quad $5$  & $8.91737$ &  $8.98792$ \\
\hline \quad $9 $  &\quad $3$  & $11.47069$ &  $11.51754$ & \quad$9$ &\quad $5$  & $11.43452 $&  $11.51754$ \\
\hline \quad $9 $  &\quad $7$  & $11.43908$ & $11.51754$ & ~\,  $11$&\quad $3$  & $14.00732 $&  $14.05335$ \\
\hline \hline
\quad $n$ & \quad $t$ &  ~$E(P_{n}^t)$ &  ~$E(C_{n})$ &~~  $n$& \quad $t$ & ~$E(P_{n}^t)$ & ~$E(C_{n})$ \\
\hline \quad $7 $  &\quad $6$  & $8.72057 $ &  $8.98792$ & \quad$9$ &\quad   $6$  & $11.38465 $ &  $11.51754$ \\
\hline ~\,   $10 $  &\quad $6$  & $12.93214 $ &  $12.94427 $ & ~\,  $11$ &\quad $6$  & $13.98336 $&  $14.05335$ \\
\hline ~\,   $13 $  &\quad $6$  & $16.55965 $ & $16.59246 $ & ~\,  $15$&\quad  $6$  & $19.12546 $&  $19.13354 $ \\
\hline
\end{tabular}
\end{center}
\caption{Values of $E(P_n^t)$ and $E(C_n)$ for $n=7,9,11,13,15$ and
some $t$.} \label{HLS20101114table3}
\end{table}

\vskip 2cm

\noindent {\bf Acknowledgement.} The authors are very grateful to
the referees for their helpful comments and suggestions, which
helped to improve the original manuscript.

\end{document}